\documentclass[11pt]{amsart} 
\usepackage{amsmath,amssymb,amsthm,mathrsfs}
\newtheorem{theorem}{Theorem}
\newtheorem{proposition}[theorem]{Proposition}
\newtheorem{corollary}[theorem]{Corollary}
\newtheorem{lemma}[theorem]{Lemma}

\begin{document}

\title{Embedded self-similar shrinkers of genus $0$}
\author{Simon Brendle}
\begin{abstract}
We confirm a well-known conjecture that the round sphere is the only compact, embedded self-similar shrinking solution to the mean curvature flow with genus $0$. More generally, we show that the only properly embedded self-similar shrinkers in $\mathbb{R}^3$ with vanishing intersection form are the sphere, the cylinder, and the plane. This answers two questions posed by T.~Ilmanen.
\end{abstract}
\address{Department of Mathematics \\ Columbia University \\ New York, NY 10027}
\dedicatory{Dedicated to Professor Leon Simon on the occasion of his seventieth birthday}
\maketitle

\section{Introduction}

This paper is concerned with self-similar shrinking solutions to the mean curvature flow in $\mathbb{R}^3$. A surface $M \subset \mathbb{R}^3$ is called a self-similar shrinker if it satisfies the equation $H = \frac{1}{2} \, \langle x,\nu \rangle$, where $\nu$ and $H$ denote the unit normal vector and the mean curvature, respectively. This condition guarantees that the surface $M$ moves by homotheties when evolved by the mean curvature flow.

The classification of self-similar solutions to geometric flows is a central problem with important implications for the analysis of singularities. Indeed, Huisken's montonicity formula \cite{Huisken1} implies that any tangent flow to a compact solution of mean curvature flow is a self-similar shrinker (see also \cite{Colding-Minicozzi-Pedersen} and \cite{Ecker}). The simplest example of a compact self-similar shrinker in $\mathbb{R}^3$ is the round sphere of radius $2$ centered at the origin. G.~Drugan \cite{Drugan} has recently constructed an example of a self-similar shrinker of genus $0$ which is immersed but fails to be embedded. Angenent \cite{Angenent} has constructed an example of an embedded self-similar shrinker of genus $1$. Moreover, N.~Kapouleas, S.~Kleene, and N.M.~M\o ller \cite{Kapouleas-Kleene-Moller} have constructed new examples of non-compact self-similar shrinkers using gluing techniques. These examples are embedded and have high genus.

A well-known conjecture asserts that the round sphere of radius $2$ should be the only embedded self-similar shrinker of genus $0$. Our main result confirms this conjecture:

\begin{theorem}
\label{main.theorem}
Let $M$ be a compact, embedded self-similar shrinker in $\mathbb{R}^3$ of genus $0$. Then $M$ is a round sphere.
\end{theorem}

In view of the examples constructed by Drugan and Angenent, the assumptions that $M$ is embedded and has genus $0$ are both necessary. In that respect, Theorem \ref{main.theorem} shares some common features with Lawson's Conjecture on embedded minimal tori in $S^3$ (cf. \cite{Lawson}). This conjecture was recently confirmed in \cite{Brendle3}; see \cite{Brendle4} for a survey.

In the noncompact case, our arguments imply the following: 

\begin{theorem}
\label{main.theorem.noncompact}
Suppose that $M$ is a properly embedded self-similar shrinker in $\mathbb{R}^3$ with the property that any two loops in $M$ have vanishing intersection number mod $2$. Then $M$ is a round sphere or a cylinder or a plane.
\end{theorem}

Theorem \ref{main.theorem.noncompact} confirms two conjectures of T.~Ilmanen, the Wiggly Plane Conjecture and the Planar Domain Conjecture (cf. \cite{Ilmanen}). We note that the topological assumption in Theorem \ref{main.theorem.noncompact} is equivalent to the condition that $M$ is homeomorphic to an open subset of $S^2$; this follows e.g. from the simple exhaustion theorem in Section 4 in \cite{Ferrer-Martin-Meeks}.

We next discuss some related results. In 1990, G.~Huisken \cite{Huisken1} proved that the round sphere is the only compact self-similar shrinking solution with positive mean curvature. Using a similar argument, Huisken was able to show that a noncompact self-similar shrinker which has bounded curvature and positive mean curvature must be a cylinder (cf. \cite{Huisken2}). Moreover, K.~Ecker and G.~Huisken proved that a self-similar shrinker which can be written as an entire graph must be a plane (cf. \cite{Ecker-Huisken}, p.~471). In a remarkable recent work, T.~Colding and W.~Minicozzi \cite{Colding-Minicozzi} proved that a self-similar shrinker which is a stable critical point of a certain entropy functional must be a sphere or a cylinder or a plane. Furthermore, T.~Colding, T.~Ilmanen, W.~Minicozzi, and B.~White recently showed that the round sphere has smallest entropy among all compact self-similar shrinkers (see \cite{Colding-Ilmanen-Minicozzi-White}).  We note that L.~Wang \cite{Wang1} has obtained a classification of self-similar shrinkers which are asymptotic to cones. X.~Wang \cite{Wang2} has proved a uniqueness result for convex translating solutions to the mean curvature flow which can be expressed as graphs over $\mathbb{R}^3$. Furthermore, we recently obtained a classification of steady gradient Ricci solitons in dimension $3$ and $4$ under a noncollapsing assumption (cf. \cite{Brendle1},\cite{Brendle2}).

We now sketch the main ideas involved in the proof of Theorem \ref{main.theorem}. Suppose that $M$ is a compact, embedded self-similar shrinker in $\mathbb{R}^3$ of genus $0$. In the first step, we show that, for any plane $P \subset \mathbb{R}^3$ which passes through the origin, the intersection  $M \cap P$ consists of a single Jordan curve which is piecewise $C^1$. This argument is inspired in part by the two-piece property for embedded minimal surfaces in $S^3$ (cf. Ros \cite{Ros}). We next prove that $M$ is star-shaped. Indeed, if $\langle \bar{x},\nu(\bar{x}) \rangle = 0$ for some point $\bar{x} \in M$, then we consider the tangent plane $P$ to $M$ at $\bar{x}$. Clearly, $P$ passes through the origin, so the intersection $M \cap P$ consists of a single Jordan curve. On the other hand, $M \cap P$ contains at least two arcs which intersect transversally at $\bar{x}$. This gives a contradiction. Having established that $M$ is star-shaped, it follows that the mean curvature of $M$ does not change sign. Huisken's theorem then implies that $M$ is a round sphere, thereby completing the proof of Theorem \ref{main.theorem}. 

The proof of Theorem \ref{main.theorem.noncompact} uses similar techniques; this is discussed in Section \ref{noncompact}.

The author is grateful to Otis Chodosh and Brian White for comments on an earlier version of this paper. This work was supported in part by the National Science Foundation under grant DMS-1201924.

\section{The key estimate}

\label{aux}

We begin by collecting some basic identities for self-similar shrinkers in $\mathbb{R}^3$. 

\begin{proposition}
\label{vector.field}
Let $\Sigma$ be a self-similar shrinker in $\mathbb{R}^3$. Moreover, suppose that $\Xi$ be a smooth vector field on $\mathbb{R}^3$, and let $\xi$ denote the projection of $\Xi$ to the tangent plane of $\Sigma$. Then 
\[\text{\rm div}_\Sigma \xi - \frac{1}{2} \, \langle x,\xi \rangle = \sum_{i=1}^2 \langle \bar{D}_{e_i} \Xi,e_i \rangle - \frac{1}{2} \, \langle x,\Xi \rangle.\] 
Here, $\bar{D}$ denotes the Levi-Civita connection on the ambient space $\mathbb{R}^3$, and $\{e_1,e_2\}$ is a local orthonormal frame on $\Sigma$. 
\end{proposition} 

\textbf{Proof.} 
Since $\Sigma$ is a self-similar shrinker, we have $H = \frac{1}{2} \, \langle x,\nu \rangle$. This implies 
\begin{align*} 
\text{\rm div}_\Sigma \xi - \frac{1}{2} \, \langle x,\xi \rangle 
&= \sum_{i=1}^2 \langle \bar{D}_{e_i} \Xi,e_i \rangle - H \, \langle \Xi,\nu \rangle - \frac{1}{2} \, \langle x,\xi \rangle \\ 
&= \sum_{i=1}^2 \langle \bar{D}_{e_i} \Xi,e_i \rangle - \frac{1}{2} \, \langle x,\nu \rangle \, \langle \Xi,\nu \rangle - \frac{1}{2} \, \sum_{i=1}^2 \langle x,e_i \rangle \, \langle \Xi,e_i \rangle \\ 
&= \sum_{i=1}^2 \langle \bar{D}_{e_i} \Xi,e_i \rangle - \frac{1}{2} \, \langle x,\Xi \rangle. 
\end{align*} 
This proves the assertion. \\

\begin{corollary}
\label{scalar.functions} 
Let $\Sigma$ be a self-similar shrinker in $\mathbb{R}^3$. Suppose that $F: \mathbb{R}^3 \to \mathbb{R}$ is a smooth function, and let $f: \Sigma \to \mathbb{R}$ denote the restriction of $F$ to $\Sigma$. Then 
\[\Delta_\Sigma f - \frac{1}{2} \, \langle x,\nabla^\Sigma f \rangle = \sum_{i=1}^2 (\bar{D}^2 F)(e_i,e_i) - \frac{1}{2} \, \langle x,\bar{\nabla} F \rangle.\] 
Here, $\bar{\nabla} F$ and $\bar{D}^2 F$ denote gradient and Hessian of $F$ with respect to the Euclidean metric, and $\{e_1,e_2\}$ is a local orthonormal frame on $\Sigma$. 
\end{corollary}

\textbf{Proof.} Apply Proposition \ref{vector.field} to the gradient vector field $\Xi = \bar{\nabla} F$. \\

It is well-known that self-similar shrinkers can be characterized as critical points of a functional. More precisely, $\Sigma$ is a self-similar shrinker if and only if $\Sigma$ is a critical point of the functional 
\[\mathscr{F}(\Sigma) = \int_\Sigma e^{-\frac{|x|^2}{4}}.\] 
Following Colding and Minicozzi, we define a differential operator $L$ on $\Sigma$ by 
\[Lf = \Delta_\Sigma f + |A|^2 \, f + \frac{1}{2} \, f - \frac{1}{2} \, \langle x,\nabla^\Sigma f \rangle\] 
(cf. \cite{Colding-Minicozzi}, equation (4.13)). The second variation of $\mathscr{F}$ is given by 
\[-\int_\Sigma e^{-\frac{|x|^2}{4}} \, f \, Lf = \int_\Sigma e^{-\frac{|x|^2}{4}} \, \Big ( |\nabla^\Sigma f|^2 - |A|^2 \, f^2 - \frac{1}{2} \, f^2 \Big ),\] 
where $f: \bar{\Sigma} \to \mathbb{R}$ is a test function which has compact support and vanishes along the boundary of $\Sigma$ (see \cite{Colding-Minicozzi}, Theorem 4.14).

We next consider a self-similar shrinker whose boundary is contained in a plane. In this case, we can use the height function as a test function in the stability inequality. This leads to the following result: 

\begin{proposition}
\label{key.estimate}
Let $\Sigma$ be a smooth surface in $\mathbb{R}^3$ with boundary $\partial \Sigma = \Gamma$, and let $k \geq 4$. Suppose that $H = \frac{1}{2} \, \langle x,\nu \rangle$ on $\Sigma \cap \{|x| \leq k\}$. Moreover, suppose that the stability inequality 
\[0 \leq -\int_\Sigma e^{-\frac{|x|^2}{4}} \, f \, Lf\] 
holds for every smooth function $f: \bar{\Sigma} \to \mathbb{R}$ which vanishes on the set $\Gamma \cup (\bar{\Sigma} \cap \{|x| \geq k\})$. Finally, we assume that $\Gamma \cap \{|x| \leq k\} \subset \{x \in \mathbb{R}^3: \langle a,x \rangle = 0\}$ for some unit vector $a \in \mathbb{R}^3$. Then 
\[\int_{\Sigma \cap \{|x| \leq \sqrt{k}\}} |A|^2 \, e^{-\frac{|x|^2}{4}} \, \langle a,x \rangle^2 \leq \frac{C}{\log k} \int_{\Sigma \cap \{\sqrt{k} \leq |x| \leq k\}} e^{-\frac{|x|^2}{4}},\] 
where $C$ is a positive constant independent of $k$.
\end{proposition}

\textbf{Proof.} 
Let us fix a smooth cutoff function $\eta: (-\infty,\infty) \to [0,1]$ satisfying $\eta = 1$ on $(-\infty,\frac{1}{2}]$, $\eta = 0$ on $[1,\infty)$, and $\eta' \leq 0$ on $(-\infty,\infty)$. We define a smooth function $F: \mathbb{R}^3 \to \mathbb{R}$ by 
\[F(x) = \langle a,x \rangle \, \eta \Big ( \frac{\log |x|}{\log k} \Big ).\] 
Note that  
\[\langle x,\bar{\nabla} F \rangle = F + \frac{1}{\log k} \, \langle a,x \rangle \, \eta' \Big ( \frac{\log |x|}{\log k} \Big ).\] 
Moreover, we have  
\[|\bar{D}^2 F| \leq \frac{C}{\log k} \, \frac{1}{|x|} \, 1_{\{\sqrt{k} \leq |x| \leq k\}},\] 
where $C$ is a positive constant independent of $k$. This implies 
\[|F| \, |\bar{D}^2 F| \leq \frac{C}{\log k} \, 1_{\{\sqrt{k} \leq |x| \leq k\}}.\] 
Let $f: \bar{\Sigma} \to \mathbb{R}$ denote the restriction of $F$ to $\bar{\Sigma}$. Using Corollary \ref{scalar.functions}, we obtain 
\begin{align*} 
&-f \, \Big ( \Delta_\Sigma f - \frac{1}{2} \, \langle x,\nabla^\Sigma f \rangle \Big ) \\ 
&= -F \, \bigg ( \sum_{i=1}^2 (\bar{D}^2 F)(e_i,e_i) - \frac{1}{2} \, \langle x,\bar{\nabla} F \rangle \bigg ) \\ 
&= -F \sum_{i=1}^2 (\bar{D}^2 F)(e_i,e_i) + \frac{1}{2} \, F^2 + \frac{1}{2\log k} \, \langle a,x \rangle^2 \, \eta \Big ( \frac{\log |x|}{\log k} \Big ) \, \eta' \Big ( \frac{\log |x|}{\log k} \Big ) \\ 
&\leq \frac{C}{\log k} \, 1_{\{\sqrt{k} \leq |x| \leq k\}} + \frac{1}{2} \, f^2. 
\end{align*} 
In the last step, we have used the inequality $\eta' \leq 0$. Consequently, 
\[-f \, Lf \leq \frac{C}{\log k} \, 1_{\{\sqrt{k} \leq |x| \leq k\}} - |A|^2 \, f^2.\] 
Note that $f$ vanishes on the set $\Gamma \cup (\bar{\Sigma} \cap \{|x| \geq k\})$. Using $f$ as a test function in the stability inequality gives 
\begin{align*} 
0 &\leq -\int_\Sigma e^{-\frac{|x|^2}{4}} \, f \, Lf \\ 
&\leq \frac{C}{\log k} \int_{\Sigma \cap \{\sqrt{k} \leq |x| \leq k\}} e^{-\frac{|x|^2}{4}} - \int_\Sigma |A|^2 \, e^{-\frac{|x|^2}{4}} \, f^2 \\ 
&\leq \frac{C}{\log k} \int_{\Sigma \cap \{\sqrt{k} \leq |x| \leq k\}} e^{-\frac{|x|^2}{4}} - \int_{\Sigma \cap \{|x| \leq \sqrt{k}\}} |A|^2 \, e^{-\frac{|x|^2}{4}} \, \langle a,x \rangle^2. 
\end{align*} 
This proves the assertion. 

\section{Proof of Theorem \ref{main.theorem}}

\label{proof.of.main.thm}

In this section, we describe the proof of Theorem \ref{main.theorem}. Let $M$ be a compact, embedded self-similar shrinker in $\mathbb{R}^3$ of genus $0$. Moreover, suppose that $M$ is not a round sphere. By Theorem 4.1 in \cite{Huisken1}, the mean curvature $H$ must change sign. In particular, we can find a point $\bar{x} \in M$ such that $H(\bar{x}) = 0$. Using the shrinker equation, we conclude that $\langle \bar{x},\nu(\bar{x}) \rangle = 0$. For abbreviation, let $a := \nu(\bar{x})$ and $Z := \{x \in M: \langle a,x \rangle = 0\}$. Clearly, $\bar{x} \in Z$. The structure of the set $Z$ is described in the following lemma.

\begin{lemma}
\label{nodal.set}
The set $Z = \{x \in M: \langle a,x \rangle = 0\}$ is a union of finitely many $C^1$-arcs which meet at isolated points. More precisely, for each point $x_0 \in Z$, there exists an open neighborhood $U \subset M$ of $x_0$ such that $Z \cap U$ is a union of $m$ $C^1$-arcs which intersect transversally at $x_0$. Here, $m$ can be characterized as the order of vanishing of the function $x \mapsto \langle a,x \rangle$ at $x_0$.
\end{lemma}

\textbf{Proof.} 
The set $Z$ can be viewed as the nodal set of a solution of an elliptic equation. Indeed, it follows from Corollary \ref{scalar.functions} that the function $f(x) := \langle a,x \rangle$ satisfies the equation 
\[\Delta_M f - \frac{1}{2} \, \langle x,\nabla^M f \rangle = -\frac{1}{2} \, f\] 
(see also \cite{Colding-Minicozzi}, Lemma 3.20). This identity can be rewritten as  
\[\Delta_M (e^{-\frac{|x|^2}{8}} \, f) = h \, e^{-\frac{|x|^2}{8}} \, f,\] 
where $h := e^{\frac{|x|^2}{8}} \, \Delta_M(e^{-\frac{|x|^2}{8}}) - \frac{1}{2}$. If we apply Lemma 2.4 and Theorem 2.5 in \cite{Cheng} to the function $e^{-\frac{|x|^2}{8}} \, f$, the assertion follows. \\

We now continue with the proof of Theorem \ref{main.theorem}. In view of our choice of $\bar{x}$ and $a$, the function $x \mapsto \langle a,x \rangle$ vanishes to order $m \geq 2$ at the point $\bar{x}$. Consequently, there exists an open neighborhood $U \subset M$ of $\bar{x}$ such that $Z \cap U$ is a union of at least two $C^1$-arcs which intersect transversally at $x_0$. In particular, $Z$ cannot be a Jordan curve. Hence, we can find a closed Jordan curve $\Gamma$ with the property that $\Gamma$ is piecewise $C^1$ and $\Gamma \subsetneq Z$. Since $M$ has genus $0$, $\Gamma$ bounds a disk in $M$.

The complement $\mathbb{R}^3 \setminus M$ has two connected components which we denote by $\Omega$ and $\tilde{\Omega}$. To fix notation, let us assume that $\Omega$ is unbounded and $\tilde{\Omega}$ is bounded.

\begin{proposition} 
\label{barrier.1}
There exists a smooth surface $\Sigma \subset \Omega$ such that $\bar{\Sigma} \setminus \Sigma = \Gamma$ and $|A|^2 = \langle x,\nu \rangle = 0$ at each point on $\Sigma$.
\end{proposition}

\textbf{Proof.} 
For $k$ sufficiently large, we denote by $\mathscr{C}_k$ the set of all embedded disks $S \subset \bar{\Omega} \cap \{|x| \leq 2k\}$ with the property that $\partial S = \Gamma$. The fact that $\Gamma$ bounds a disk in $M$ implies that $\mathscr{C}_k$ is non-empty if $k$ is sufficiently large. Moreover, we choose a smooth cutoff function $\psi_k: [0,\infty) \to [0,1]$ satisfying $\psi_k=0$ on $[0,k]$ and $\psi_k'(2k) > k$. We now consider the functional 
\[\mathscr{F}_k(S) = \int_S e^{-\frac{|x|^2}{4} + \psi_k(|x|)}\] 
for $S \in \mathscr{C}_k$. We can interpret $\mathscr{F}_k$ as the area functional for the conformal metric $e^{-\frac{|x|^2}{4}+\psi_k(|x|)} \, \delta_{ij}$. For $k$ sufficiently large, the region $\bar{\Omega} \cap \{|x| \leq 2k\}$ is a mean convex domain with respect to this conformal metric. Therefore, 
general results from \cite{Meeks-Yau} guarantee that there exists a smooth embedded surface $\Sigma_k \in \mathscr{C}_k$ which minimizes the functional $\mathscr{F}_k$. Since $\Sigma_k$ is a global minimizer for the functional $\mathscr{F}_k$, it is easy to see that 
\[\sup_k \mathscr{F}_k(\Sigma_k) < \infty.\] 
This implies 
\[\sup_k \int_{\Sigma_k} e^{-\frac{|x|^2}{4}} < \infty.\] 
Using the first variation formula, we deduce that $H = \frac{1}{2} \, \langle x,\nu \rangle$ on $\Sigma_k \cap \{|x| \leq k\}$. Finally, the stability inequality implies that 
\[0 \leq -\int_{\Sigma_k} e^{-\frac{|x|^2}{4}} \, f \, Lf\] 
for every smooth function $f: \bar{\Sigma}_k \to \mathbb{R}$ which vanishes on the set $\Gamma \cup (\bar{\Sigma}_k \cap \{|x| \geq k\})$. Using Proposition \ref{key.estimate}, we obtain
\begin{align} 
\label{almost.totally.geodesic}
&\limsup_{k \to \infty} \int_{\Sigma_k \cap \{|x| \leq \sqrt{k}\}} |A|^2 \, e^{-\frac{|x|^2}{4}} \, \langle a,x \rangle^2 \notag \\ 
&\leq \limsup_{k \to \infty} \frac{C}{\log k} \int_{\Sigma_k \cap \{\sqrt{k} \leq |x| \leq k\}} e^{-\frac{|x|^2}{4}} = 0. 
\end{align} 
Finally, it follows from Theorem 3 in \cite{Schoen-Simon} that 
\[\limsup_{k \to \infty} \sup_{\Sigma_k \cap W} |A|^2 < \infty\] 
for every compact set $W \subset \mathbb{R}^3 \setminus \Gamma$. Hence, after passing to a subsequence if necessary, the surfaces $\Sigma_k$ converge in $C_{loc}^\infty(\mathbb{R}^3 \setminus \Gamma)$ to a smooth surface $\Sigma \subset \mathbb{R}^3 \setminus \Gamma$ which satisfies the shrinker equation $H = \frac{1}{2} \, \langle x,\nu \rangle$. Using (\ref{almost.totally.geodesic}), we conclude that $\Sigma$ is totally geodesic. In particular, $\langle x,\nu \rangle = 0$ at each point on $\Sigma$. Moreover, it is easy to see that $\Sigma \subset \bar{\Omega}$. Since $M$ is not totally geodesic, the strict maximum principle implies that $\Sigma$ cannot touch $M$. Consequently, $\Sigma \subset \Omega$. 

We next show that $\Gamma \subset \bar{\Sigma}$. If $\Gamma \setminus \bar{\Sigma} \neq \emptyset$, we can construct a one-form $\alpha$ on $\mathbb{R}^3$ such that $\alpha$ has compact support, $\alpha = 0$ in an open neighborhood of $\bar{\Sigma}$, $d\alpha = 0$ in an open neighborhood of $\Gamma$, and $\int_\Gamma \alpha \neq 0$. This implies $\int_{\Sigma_k} d\alpha = \int_\Gamma \alpha \neq 0$ for each $k$, and $\int_{\Sigma_k} d\alpha \to 0$ as $k \to \infty$. This is a contradiction. Thus, $\Gamma \subset \bar{\Sigma}$. Since $\bar{\Sigma} \setminus \Sigma \subset \Gamma$, we conclude that $\bar{\Sigma} \setminus \Sigma = \Gamma$. This completes the proof of Proposition \ref{barrier.1}. \\


\begin{proposition}
\label{barrier.2}
There exists a smooth surface $\tilde{\Sigma} \subset \tilde{\Omega}$ such that $\bar{\tilde{\Sigma}} \setminus \tilde{\Sigma} = \Gamma$ and $|A|^2 = \langle x,\nu \rangle = 0$ at each point on $\tilde{\Sigma}$.
\end{proposition}

\textbf{Proof.} 
We consider the set $\tilde{\mathscr{C}}$ of all embedded disks $S \subset \bar{\tilde{\Omega}}$ with boundary $\partial S = \Gamma$. As above, $\tilde{\mathscr{C}}$ is non-empty since $\Gamma$ bounds a disk in $M$. We now consider the functional $\mathscr{F}$ defined in Section \ref{aux}. The functional $\mathscr{F}$ can be viewed as the area functional for the conformal metric $e^{-\frac{|x|^2}{4}} \, \delta_{ij}$. Clearly, $\bar{\tilde{\Omega}}$ is a mean convex domain with respect to this conformal metric. Using results in \cite{Meeks-Yau}, we can find a smooth embedded surface $\tilde{\Sigma} \in \tilde{\mathscr{C}}$ which minimizes the functional $\mathscr{F}$. The first variation formula implies that the surface $\tilde{\Sigma}$ satisfies $H = \frac{1}{2} \, \langle x,\nu \rangle$. Moreover, the stability inequality gives 
\[0 \leq -\int_{\tilde{\Sigma}} e^{-\frac{|x|^2}{4}} \, f \, Lf\] 
for every smooth function $f: \bar{\tilde{\Sigma}} \to \mathbb{R}$ which vanishes on the boundary $\Gamma$. Using Proposition \ref{key.estimate} with $k$ sufficiently large, we obtain 
\[\int_{\tilde{\Sigma}} |A|^2 \, e^{-\frac{|x|^2}{4}} \, \langle a,x \rangle^2 = 0.\] 
Consequently, $\tilde{\Sigma}$ is totally geodesic. This implies $\langle x,\nu \rangle = 0$ at each point on $\tilde{\Sigma}$. Finally, we clearly have $\tilde{\Sigma} \subset \bar{\tilde{\Omega}}$. Since $M$ is not totally geodesic, the strict maximum principle implies that $\tilde{\Sigma}$ cannot touch $M$. Therefore, $\tilde{\Sigma} \subset \tilde{\Omega}$, as claimed. \\

\begin{proposition}
\label{normal.vector.of.Sigma.and.tilde.Sigma}
The unit normal vectors to $\Sigma$ and $\tilde{\Sigma}$ are parallel to $a$ at all points. 
\end{proposition} 

\textbf{Proof.} 
Suppose that there exists a point $x \in \Sigma$ such that $\nu(x) = b$, where $a$ and $b$ are linearly independent. Let us define 
\[\Sigma' = \{x \in \Sigma: \nu(x)=b\} \neq \emptyset.\] 
By Proposition \ref{barrier.1}, $\Sigma'$ is a subset of $\{x \in \mathbb{R}^3: \langle b,x \rangle = 0\} \setminus \Gamma$. Moreover, Proposition \ref{barrier.1} implies that $\Sigma'$, viewed as a subset of $\{x \in \mathbb{R}^3: \langle b,x \rangle = 0\} \setminus \Gamma$, is open and closed. On the other hand, we have 
\[\{x \in \mathbb{R}^3: \langle b,x \rangle = 0\} \cap \Gamma \subset \{x \in \mathbb{R}^3: \langle a,x \rangle = \langle b,x \rangle = 0\} =: L.\] 
Hence, the closure of $\Sigma'$ is either an entire plane or a halfplane with boundary $L$. In the latter case, $\Gamma$ contains a straight line, but this is impossible since $M$ is compact. Consequently, the closure of $\Sigma'$ is the entire plane $\{x \in \mathbb{R}^3: \langle b,x \rangle = 0\}$. Since $\Sigma' \subset \Omega$, it follows that the plane $\{x \in \mathbb{R}^3: \langle b,x \rangle = 0\}$ is contained in $\bar{\Omega}$, and $M$ lies on one side of this plane. This contradicts the fact that $\int_M e^{-\frac{|x|^2}{4}} \, \langle b,x \rangle = 0$. Consequently, the normal vector to $\Sigma$ is parallel to $a$ at each point on $\Sigma$. An analogous argument shows that the normal vector to $\tilde{\Sigma}$ is parallel to $a$ at each point on $\tilde{\Sigma}$. This completes the proof of Proposition \ref{normal.vector.of.Sigma.and.tilde.Sigma}. \\

Combining Propositions \ref{barrier.1}, \ref{barrier.2}, and \ref{normal.vector.of.Sigma.and.tilde.Sigma}, we conclude that the surfaces $\Sigma$ and $\tilde{\Sigma}$ are contained in the plane $\{x \in \mathbb{R}^3: \langle a,x \rangle = 0\}$. Moreover, we have $\Sigma \subset \Omega$ and $\tilde{\Sigma} \subset \tilde{\Omega}$; in particular, $\Sigma$ and $\tilde{\Sigma}$ are disjoint. Finally, $\Sigma$ and $\tilde{\Sigma}$ have the same boundary $\Gamma$. Therefore, the union $\Sigma \cup \tilde{\Sigma} \cup \Gamma$, viewed as a subset of $\{x \in \mathbb{R}^3: \langle a,x \rangle = 0\}$, is open and closed. This implies 
\[\{x \in \mathbb{R}^3: \langle a,x \rangle = 0\} = \Sigma \cup \tilde{\Sigma} \cup \Gamma \subset \Omega \cup \tilde{\Omega} \cup \Gamma = (\mathbb{R}^3 \setminus M) \cup \Gamma.\] 
Consequently,  
\[\{x \in M: \langle a,x \rangle = 0\} \subset \Gamma.\] 
In other words, the set $Z$ coincides with $\Gamma$. This contradicts our choice of $\Gamma$. This completes the proof of Theorem \ref{main.theorem}.

\section{Proof of Theorem \ref{main.theorem.noncompact}} 

\label{noncompact}

In this final section, we discuss the proof of Theorem \ref{main.theorem.noncompact}. Throughout this section, we assume that $M$ is a properly embedded self-similar shrinker in $\mathbb{R}^3$. We first recall a well known result, which is an immediate consequence of Brakke's local area bound (see \cite{Brakke} or \cite{Ecker}, Proposition 4.9): 

\begin{proposition} 
\label{area.growth}
For $k$ large, the area of $M \cap \{|x| \leq k\}$ is at most $O(k^2)$. 
\end{proposition}

\textbf{Proof.} 
We sketch the proof for the convenience of the reader. By assumption, the surfaces $M_t = \sqrt{-t} \, M$ form a solution of mean curvature flow. Applying Proposition 4.9 in \cite{Ecker} (with $t_0=0$ and $\rho=4$) gives 
\[\text{\rm area}(M_t \cap \{|x| \leq 2\}) \leq 8 \, \text{\rm area}(M_{-1} \cap \{|x| \leq 4\})\] 
for all $t \in [-1,0)$. This implies  
\[\text{\rm area}(M \cap \{|x| \leq 2k\}) \leq 8k^2 \, \text{\rm area}(M \cap \{|x| \leq 4\})\] 
for $k \geq 1$. From this, the assertion follows. \\

We will also need the following result, which is a special case of a much stronger theorem due to Brian White:

\begin{theorem}[B.~White \cite{White}]
\label{brian.white.theorem}
Suppose that $M$ contains the line $\{x \in \mathbb{R}^3: x_1=x_2=0\}$, and $M$ is disjoint from the halfplane $\{x \in \mathbb{R}^3: x_1 < 0, \, x_2=0\}$. Then $M$ is a plane.
\end{theorem}

\textbf{Proof.} 
We again sketch an argument for the convenience of the reader. We define a vector field $\Xi$ on $\mathbb{R}^3$ by $\Xi(x_1,x_2,x_3) = (-x_2,x_1,0)$, and let $\xi$ denote the projection of $\Xi$ to the tangent plane of $M$. By assumption, $M$ is disjoint from the halfplane $\{x \in \mathbb{R}^3: x_1<0, \, x_2=0\}$. Hence, every point $x \in M \setminus \{x \in \mathbb{R}^3: x_1=x_2=0\}$ can be uniquely written in the form $x = (\sqrt{x_1^2+x_2^2} \, \cos \theta,\sqrt{x_1^2+x_2^2} \, \sin \theta,x_3)$ for some $\theta \in (-\pi,\pi)$. This defines a smooth function $\theta: M \setminus \{x \in \mathbb{R}^3: x_1=x_2=0\} \to (-\pi,\pi)$ satisfying $(x_1^2+x_2^2) \, \nabla^M \theta = \xi$. 

Using Proposition \ref{vector.field}, we obtain 
\[\text{\rm div}_M \xi - \frac{1}{2} \, \langle x,\xi \rangle = \sum_{i=1}^2 \langle \bar{D}_{e_i} \Xi,e_i \rangle - \frac{1}{2} \, \langle x,\Xi \rangle = 0.\] 
In the last step, we have used that $\Xi$ is a Killing vector field in ambient space. This implies 
\begin{align*} 
\text{\rm div}_M(e^{-\frac{|x|^2}{4}} \, \theta \, \xi) 
&= e^{-\frac{|x|^2}{4}} \, \theta \, \Big ( \text{\rm div}_M \xi - \frac{1}{2} \, \langle x,\xi \rangle \Big ) + e^{-\frac{|x|^2}{4}} \, \langle \nabla^M \theta,\xi \rangle \\ 
&= e^{-\frac{|x|^2}{4}} \, (x_1^2+x_2^2) \, |\nabla^M \theta|^2 
\end{align*} 
on $M \setminus \{x \in \mathbb{R}^3: x_1=x_2=0\}$. Integrating over $M \setminus \{x \in \mathbb{R}^3: x_1=x_2=0\}$ gives 
\begin{align*} 
0 &= \int_{M \setminus \{x \in \mathbb{R}^3: x_1=x_2=0\}} \text{\rm div}_M(e^{-\frac{|x|^2}{4}} \, \theta \, \xi) \\ 
&= \int_{M \setminus \{x \in \mathbb{R}^3: x_1=x_2=0\}} e^{-\frac{|x|^2}{4}} \, (x_1^2+x_2^2) \, |\nabla^M \theta|^2. 
\end{align*}
Consequently, $\nabla^M \theta = 0$ at each point $x \in M \setminus \{x \in \mathbb{R}^3: x_1=x_2=0\}$. Thus, $M$ is a plane. This completes the proof of Theorem \ref{brian.white.theorem}. \\

We now continue with the proof of Theorem \ref{main.theorem.noncompact}. Let $M$ be a properly embedded self-similar shrinker in $\mathbb{R}^3$ with the property that any two loops in $M$ have vanishing intersection number mod $2$. Moreover, we assume that $M$ is neither a round sphere, nor a cylinder, nor a plane. By Proposition \ref{area.growth}, $M$ has polynomial area growth. By a theorem of Colding and Minicozzi, the mean curvature must change sign (see \cite{Colding-Minicozzi}, Theorem 10.1). In particular, we can find a point $\bar{x} \in M$ such that $H(\bar{x}) = 0$, and hence $\langle \bar{x},\nu(\bar{x}) \rangle = 0$. As above, we put $a := \nu(\bar{x})$. We now consider two cases: 

\textit{Case 1:} Suppose that the sets $\{x \in M: \langle a,x \rangle > 0\}$ and $\{x \in M: \langle a,x \rangle < 0\}$ are both connected. In this case, we can construct two loops with the property that the first loop is contained in $\{x \in M: \langle a,x \rangle > 0\} \cup \{\bar{x}\}$, the second loop is contained in $\{x \in M: \langle a,x \rangle < 0\} \cup \{\bar{x}\}$, and the two loops intersect transversally at $\bar{x}$. This contradicts our assumption that any two loops in $M$ have vanishing intersection number mod $2$.

\textit{Case 2:} For the remainder of this section, we will assume that one of the sets $\{x \in M: \langle a,x \rangle > 0\}$ and $\{x \in M: \langle a,x \rangle < 0\}$ is not connected. Without loss of generality, we may assume that $\{x \in M: \langle a,x \rangle > 0\}$ is disconnected. Let $D$ be an arbitrary connected component of $\{x \in M: \langle a,x \rangle > 0\}$. Moreover, let $\Omega$ and $\tilde{\Omega}$ denote the connected components of $\mathbb{R}^3 \setminus M$. 

\begin{proposition} 
\label{barrier.1.sec.4}
There exists a smooth surface $\Sigma \subset \Omega$ such that $\bar{\Sigma} \setminus \Sigma = \partial D$ and $|A|^2 = \langle x,\nu \rangle = 0$ at each point on $\Sigma$.
\end{proposition} 

\textbf{Proof.} 
By Sard's Lemma, we can find a sequence of numbers $r_k \in (2k,3k)$ such that the sphere $\{|x|=r_k\}$ intersects $M$ transversally. By smoothing out the domain $\Omega \cap \{|x| < r_k\}$, we can construct an open domain $\Omega_k$ with smooth boundary such that $\Omega_k \cap \{|x| \leq 2k\} = \Omega \cap \{|x| \leq 2k\}$ and $\Omega_k \subset \{|x| \leq 3k\}$. Moreover, we can find a smooth function $\chi_k: \bar{\Omega}_k \to [0,1]$ such that $\chi_k = 0$ on the set $\bar{\Omega}_k \cap \{|x| \leq k\}$ and $\bar{\Omega}_k$ is a mean convex domain with respect to the conformal metric $e^{-\frac{|x|^2}{4}+\chi_k(x)} \, \delta_{ij}$. 

By Sard's Lemma, we can find a real number $\rho_k \in (k,2k)$ such that the sphere $\{|x|=\rho_k\}$ intersects $M$ and $\partial D$ transversally. Clearly, the curve $\Gamma_k = \partial (D \cap \{|x| < \rho_k\})$ satisfies $\Gamma_k \cap \{|x| \leq k\} \subset \partial D \subset \{x \in \mathbb{R}^3: \langle a,x \rangle = 0\}$. Let $\Sigma_k$ be a surface which minimizes the modified area functional 
\[\int_S e^{-\frac{|x|^2}{4}+\chi_k(x)}\] 
among all embedded, orientable surfaces $S \subset \bar{\Omega}_k$ with boundary $\partial S = \Gamma_k$. Clearly, 
\[\sup_k \int_{\Sigma_k} e^{-\frac{|x|^2}{4}} < \infty.\] 
Moreover, the first variation formula implies that $H = \frac{1}{2} \, \langle x,\nu \rangle$ at each point on $\Sigma_k \cap \{|x| \leq k\}$. Using the stability inequality together with Proposition \ref{key.estimate}, we conclude that 
\begin{align} 
\label{almost.totally.geodesic.sec.4}
&\limsup_{k \to \infty} \int_{\Sigma_k \cap \{|x| \leq \sqrt{k}\}} |A|^2 \, e^{-\frac{|x|^2}{4}} \, \langle a,x \rangle^2 \notag \\ 
&\leq \limsup_{k \to \infty} \frac{C}{\log k} \int_{\Sigma_k \cap \{\sqrt{k} \leq |x| \leq k\}} e^{-\frac{|x|^2}{4}} = 0. 
\end{align} 
Finally, it follows from results in \cite{Schoen-Simon} that  
\[\limsup_{k \to \infty} \sup_{\Sigma_k \cap W} |A|^2 < \infty\] 
for every compact set $W \subset \mathbb{R}^3 \setminus \partial D$. Hence, after passing to a subsequence, the surfaces $\Sigma_k$ converge in $C_{loc}^\infty(\mathbb{R}^3 \setminus \partial D)$ to a smooth surface $\Sigma \subset \mathbb{R}^3 \setminus \partial D$ which satisfies the shrinker equation $H = \frac{1}{2} \, \langle x,\nu \rangle$. Using (\ref{almost.totally.geodesic.sec.4}), we conclude that $\Sigma$ is totally geodesic. In particular, $\langle x,\nu \rangle = 0$ at each point on $\Sigma$. Moreover, it is easy to see that $\Sigma \subset \bar{\Omega}$. Since $M$ is not totally geodesic, the strict maximum principle implies that $\Sigma$ cannot touch $M$. Consequently, $\Sigma \subset \Omega$. Arguing as above, we obtain $\bar{\Sigma} \setminus \Sigma = \partial D$. \\

\begin{proposition}
\label{barrier.2.sec.4} 
There exists a smooth surface $\tilde{\Sigma} \subset \Omega$ such that $\bar{\tilde{\Sigma}} \setminus \tilde{\Sigma} = \partial D$ and $|A|^2 = \langle x,\nu \rangle = 0$ at each point on $\tilde{\Sigma}$.
\end{proposition}

\textbf{Proof.} 
Analogous to Proposition \ref{barrier.1.sec.4}. \\

\begin{proposition}
\label{normal.vector.of.Sigma.and.tilde.Sigma.sec.4}
The unit normal vectors to $\Sigma$ and $\tilde{\Sigma}$ are parallel to $a$ at all points. 
\end{proposition} 

\textbf{Proof.} 
Suppose that there exists a point on $\Sigma$ or $\tilde{\Sigma}$ where the unit normal vector is not parallel to $a$. Without loss of generality, we may assume that there exists a point $x \in \Sigma$ such that $\nu(x)=b$, where $a$ and $b$ are linearly independent. Let 
\[\Sigma' = \{x \in \Sigma: \nu(x)=b\} \neq \emptyset.\] 
As above, $\Sigma'$ is a subset of $\{x \in \mathbb{R}^3: \langle b,x \rangle = 0\} \setminus \partial D$, which is both open and closed. Since 
\[\{x \in \mathbb{R}^3: \langle b,x \rangle = 0\} \cap \partial D \subset \{x \in \mathbb{R}^3: \langle a,x \rangle = \langle b,x \rangle = 0\} =: L,\] 
it follows that the closure of $\Sigma'$ is either an entire plane or a halfplane with boundary $L$. If the closure of $\Sigma'$ is a halfplane with boundary $L$, then $M$ contains the line $L$, and $M$ is disjoint from $\Sigma' \setminus L$. Hence, it follows from Theorem \ref{brian.white.theorem} that $M$ is a plane, contrary to our assumption. Since $\Sigma' \subset \Omega$, the closure of $\Sigma'$ is the entire plane $\{x \in \mathbb{R}^3: \langle b,x \rangle = 0\}$, and $M$ lies on one side of this plane. As above, this contradicts the fact that $\int_M e^{-\frac{|x|^2}{4}} \, \langle b,x \rangle = 0$. This completes the proof. \\

It follows from Propositions \ref{barrier.1.sec.4}, \ref{barrier.2.sec.4}, and \ref{normal.vector.of.Sigma.and.tilde.Sigma.sec.4} that $\Sigma$ and $\tilde{\Sigma}$ are contained in the plane $\{x \in \mathbb{R}^3: \langle a,x \rangle = 0\}$. Moreover, $\Sigma$ and $\tilde{\Sigma}$ are disjoint, and have the same boundary $\partial D$. Putting these facts together, we conclude that 
\[\{x \in \mathbb{R}^3: \langle a,x \rangle = 0\} = \Sigma \cup \tilde{\Sigma} \cup \partial D \subset \Omega \cup \tilde{\Omega} \cup \partial D = (\mathbb{R}^3 \setminus M) \cup \partial D.\] 
Thus, $\{x \in M: \langle a,x \rangle = 0\} = \partial D$. This implies $\{x \in M: \langle a,x \rangle > 0\} = D$. In particular, the set $\{x \in M: \langle a,x \rangle > 0\}$ is connected, contrary to our assumption. This completes the proof of Theorem \ref{main.theorem.noncompact}.

\end{document}